\documentclass[10pt]{amsart}
\usepackage[cp1251]{inputenc}
\usepackage[english,russian]{babel}
\usepackage{amsmath}
\usepackage{amssymb}
\usepackage{amsfonts}
\usepackage{graphicx}

\def\mscs{51M04,~12F10} %Здесь автор определяет классификаторы AMS своей работы

\newtheorem{lemma}{Lemma}
\newtheorem{theorem}{Theorem}

\newtheorem{corollary}{Corollary}

\begin{document}
\renewcommand{\refname}{References}

\thispagestyle{empty}

\title[Area of a triangle and angle bisectors]{Area of a triangle and angle bisectors}
\author{{A.A. Buturlakin}}%
\address{Aleksandr Aleksandrovich Buturlakin  % Контактные данные всех авторов и место работы указываются только на английском языке
\newline\hphantom{iii} Sobolev Institute of Mathematics,
\newline\hphantom{iii} Koptyuga av., 4,
\newline\hphantom{iii} 630090, Novosibirsk, Russia%
\newline\hphantom{iii} Novosibirsk State University,
\newline\hphantom{iii} Pirogova str., 1,
\newline\hphantom{iii} 630090, Novosibirsk, Russia}%
\email{buturlakin@math.nsc.ru}%

\author{{S.S. Presnyakov}}%
\address{Sergei Sergeevich Presnyakov  % Контактные данные всех авторов и место работы указываются только на английском языке
\newline\hphantom{iii} Specialized Educational Scientific Center of Novosibirsk State
University,
\newline\hphantom{iii} Pirogova str., 11/1,
\newline\hphantom{iii} 630090, Novosibirsk, Russia}%
\email{sier.pr@mail.ru}%

\author{{D.O. Revin}}%
\address{Danila Olegovich Revin  % Контактные данные всех авторов и место работы указываются только на английском языке
\newline\hphantom{iii} Sobolev Institute of Mathematics,
\newline\hphantom{iii}  Koptyuga av., 4,
\newline\hphantom{iii} 630090, Novosibirsk, Russia%
\newline\hphantom{iii} Novosibirsk State University,
\newline\hphantom{iii} Pirogova str., 1,
\newline\hphantom{iii} 630090, Novosibirsk, Russia}%
\email{revin@math.nsc.ru}%

\author{{S.A. Savin}}%
\address{Semen Antonovich Savin  % Контактные данные всех авторов и место работы указываются только на английском языке
\newline\hphantom{iii} Specialized Educational Scientific Center of Novosibirsk State
University,
\newline\hphantom{iii} Pirogova str., 11/1
\newline\hphantom{iii} 630090, Novosibirsk, Russia%
\newline\hphantom{iii} The Orthodox Gymnasium in the name saint Sergius of Radonezh.
\newline\hphantom{iii} Akademicheskaya str., 3
\newline\hphantom{iii} 630090, Novosibirsk, Russia}%
\email{semen.savin.02@gmail.com}%

%\thanks{\sc Ivanov, I.I.,
%On some problems of commutative algebra}
%\thanks{\copyright \ 2015 Иванов И.И}
\thanks{\rm Funding: The reported study was funded by RFBR and BRFBR, project number  20-51-00007, by Mathematical Center in Akademgorodok under agreement No 075-2019-1675 with the Ministry of Science and Higher Education of the Russian Federation, and by the Program of Fundamental Scientific Research of the SB RAS No. I.1.1., project  number 0314-2016-0001. }
%\thanks{\it Submitted May 6, 2020, published ??.}%

 \hphantom{aaaaaaaaaaaaaaaaaaaaaaaaaaaaaaaaaaaaaaaaaaaaaaaaaaaaaaaaaaaaaaaaaaaaaaaaaaaaaaaaaaaaaaaaaaa}{\rm\small MSC\ \ \mscs }
%\semrtop \vspace{1cm}
\maketitle 
%\vspace{1mm}
     %{\flushleft\it Vol. 17, стр. ??--?? (2020) 
     %\hspace{65mm}{\rm\small UDC \udcs}
     %\newline
     	%{\rm\small DOI~10.33048/semi.2020.17.xxx}

{\small
\begin{quote}
\noindent{\sc Abstract. }  Consider a triangle $ABC$ with given lengths $l_a,l_b,l_c$ of its internal angle bisectors. We prove that in general, it is impossible to construct a square of the same area
as $ABC$ using a ruler and compass. Moreover, it is impossible to express the area of $ABC$ in radicals of $l_a,l_b,l_c$.\medskip

 \end{quote}
}

\bigskip

Consider a triangle $ABC$ on a Euclidean plane. Let $BC=a$, $AC=b$ and $AB=c$. By $m_a,\ m_b,\ m_c$ we denote the medians of $ABC$ to the sides $a$, $b$ and $c$ respectively. Similarly, $h_a,\ h_b,\
h_c$ are the altitudes of $ABC$ and $l_a,\ l_b,\ l_c$ are the angle bisectors.

Recall that according to Heron's formula, the area $S$ of $ABC$ is equal to
$$
S=\frac{1}{4}\sqrt{(a+b+c)(a+b-c)(b+c-a)(c+a-b)}=\sqrt{p(p-a)(p-b)(p-c)},
$$
where $p=(a+b+c)/2$ is the semiperimeter of the triangle. There are equivalent formulae which express $S$ in terms of the medians
$$
S=\frac{1}{3}\sqrt{\left(m_a+m_b+m_c\right)\left(m_a+m_b-m_c\right)\left(m_b+m_c-m_a\right)\left(m_c+m_a-m_b\right)},
$$
 or the altitudes
$$
\frac{1}{S}=\sqrt{\left(\frac{1}{h_a}+\frac{1}{h_b}+\frac{1}{h_c}\right)\left(\frac{1}{h_a}+\frac{1}{h_b}-\frac{1}{h_c}\right)\left(\frac{1}{h_b}+\frac{1}{h_c}-\frac{1}{h_a}\right)\left(\frac{1}{h_c}
+\frac{1}{h_a}-\frac{1}{h_b}\right)}.
$$
Here we answer the question of {\it whether it is possible to express the area of a triangle (for example, in radicals) using the lengths of its internal angle bisectors}. This question was posed in
~\cite[problem~12]{Pluecker} (published in 1830), but undoubtedly had been known before.

The following facts should also be kept in mind:

\begin{itemize}
 \item[(i)] If the lengths of corresponding angle bisectors of two triangles are equal then the triangles are congruent.

\item[(ii)] For every triplet of positive numbers $l_a$, $l_b$, $l_c$, there exists a triangle with the lengths of the angle bisectors equal to these numbers \cite{MironescuPanaitopol}. This statement
answers in the affirmative the question by A. Brocard  \cite{Brocard}.

\item[(iii)] Two previous statements yield that there is a function
$$
(l_a,l_b,l_c)\mapsto S,
$$
 which maps every triplet of positive numbers to the area of a triangle with given angle bisectors.

\item[(iv)] In general case, it is impossible to construct a triangle given the lengths of its angle bisectors using a ruler and compass (P.~Barbarin, 1896 \cite{Barbarin},
A.~Korselt, 1897 \cite{Korselt}).

\end{itemize}

In Heron's time, the question of expressing the area of a triangle in terms of the lengths of its angle bisectors could have been posed as follows: {\it Is it possible to construct a square with the
same area as a triangle with given lengths of angle bisectors using a ruler and compass?} The following theorem answers the question negatively.

\begin{theorem}\label{T1}
{\it  In general case, it is impossible to construct a square of the same area as a triangle with given lengths of angle bisectors using a ruler and compass.}
 \end{theorem}

 PROOF. Consider a triangle $ABC$ with angle bisectors $l_a=AL=1$, $l_b=BL=1/3$ and $l_c=CL=1/3$. From~(ii) it follows that such a triangle exists. Furthermore,~(i) implies that $ABC$ is an isosceles
triangle ($AB=AC$). To prove the theorem, it is sufficient to prove that the construction is impossible for this choice of lengths. It is easy to see that the area of $ABC$ is equal
to~$\mathrm{tg}(A/2)$, so a side of the required square equals to $\sqrt{\mathrm{tg}( A/2)}$.
Moreover, given a segment of length 1, the problems of constructing the following objects are equivalent:

   \medskip

 \begin{minipage}[c]{70mm}
 \begin{itemize}
   \item a segment of length $\sqrt{\mathrm{tg}( A/2)}$;
   \item a segment of length ${\mathrm{tg}( A/2)}$;
   \item an angle $ A$;
   \item a triangle $\triangle ABC$;
   \item an angle $B/2$;
   \item a segment of length ${\sin(B/2)}$.
 \end{itemize}
  \end{minipage}\!\!\!\!
   %\begin{minipage}[c]{30mm}
%   \includegraphics{manin_triangle}
%   \end{minipage}
%
   \medskip

 Therefore, the problem of constructing a square of area equal to the area of $ABC$ using a ruler and compass is equivalent to the problem of constructing an isosceles triangle with angle bisectors
of lengths $1$, $1/3$, $1/3$  and is also equivalent to the problem of constructing a half base angle of that triangle.

It is known that in a triangle with angles $\alpha,\ \beta,\ \gamma$ and
semiperimeter $p$ the length of the bisector $l$ of angle $\alpha$ is equal to

 $$
 l=2p\displaystyle\frac{\sin\displaystyle\frac{\beta}{2}\sin\displaystyle\frac{\gamma}{2}}{\cos\displaystyle\frac{\alpha}{2}\cos\displaystyle
 \frac{\beta-\gamma}{2}}.
 $$
 In our case,
$$
 1=2p\displaystyle\frac{\sin^2\displaystyle\frac{ B}{2}}{\cos\displaystyle\frac{A}{2}}=2p\displaystyle\frac{\sin^2\displaystyle\frac{B}{2}}{\sin B}=p\,{\mathrm{tg}(B/2)} $$
 and, since $ B= C=\displaystyle\frac{\pi}{2}-\displaystyle\frac{ A}{2},$
$$\displaystyle
 \frac{1}{3}=2p\displaystyle\frac{\sin\displaystyle\frac{ A}{2}\sin\displaystyle\frac{ B}{2}}{\cos\displaystyle\frac{ B}{2}\cos\displaystyle\frac{1}{2}\left(\pi-2 B-
B\right)}=2p\displaystyle\frac{\cos B\sin\displaystyle\frac{ B}{2}}{\cos\displaystyle\frac{ B}{2}\sin\displaystyle\frac{3B}{2} }=2p\displaystyle\frac{\cos B}{\sin\displaystyle\frac{3B}{2}
}{\mathrm{tg}\displaystyle\frac{ B}{2}}.
 $$
 Hence $\displaystyle\frac{1}{3}=2\displaystyle\frac{\cos B}{\sin\displaystyle\frac{3B}{2} }$ or, in other words, $\sin(3 B/2)=6\cos B$. Let $\sin(B/2)=x$. Then $\sin(3 B/2)=3x-4x^3$, $\cos B=1-2x^2$.
We see that $3x-4x^3=6(1-2x^2)$ and
 $$
 4x^3-12 x^2-3x+6=0.
 $$
 Recall the following statement from Galois theory (see, for instance, \cite[Theorem~2]{Manin} or \cite[page~199]{Postnikov}): {\em Consider segments $l_1,\dots l_n$ of rational lengths and assume
that a segment, with its length $l$ equal to a root of some irreducible polynomial of degree $n$ with rational coefficients, can be constructed from $l_1 \ldots l_n$ using a ruler and compass. Then
$n$ is a power of $2$.}
The polynomial $4x^3-12 x^2-3x+6$ is irreducible by Eisenstein's criterion~\cite{Eisenstein}. Therefore, a segment of length $\sin( B/2)$ (which is a root of that polynomial) cannot be constructed
using a ruler and compass. The theorem is proved.

\bigskip

  Note that our proof of theorem~\ref{T1} contains a proof of the fact that a triangle with bisectors of length 1, 1/3, 1/3 cannot be constructed using a  ruler and compass. This proof can also be
found in \cite[pages 224--225]{Manin} and \cite[pages~201--202]{Postnikov} and is listed here for the sake of completeness.

 Theorem ~\ref{T1} yields the statement which strengthens the result of Barbarin and Korselt.

\begin{corollary}\label{C1}
%\noindent{\bf ТЕОРЕМА.}
{\it In general case, it is impossible to construct a triangle of the same area as a triangle with given lengths of angle bisectors using a  ruler and compass.}

 \end{corollary}

 \begin{corollary}\label{C2}
%\noindent{\bf ТЕОРЕМА.}
{\it There exists no formula, to express the area of a triangle in terms of the lengths of its angle bisectors in quadratic radicals.}

 \end{corollary}

 Concerning the question of whether there exist an explicit expression of a triangle's area in radicals of the lengths of its angle bisectors, a certain confusion persists in mathematical literature.
For instance, in \cite[page~335]{DincaMawhin} it is stated that von Renthe Fink found such an expression in 1843 \cite{vRF}. However, his expression  includes the radius of the incircle $r$ in additon
to $l_a$, $l_b$, $l_c$. One of the relations he found is as follows \cite[page~274, equation 6]{vRF}:

 $$
 4a_2r^2S^2-8a_3r^3S^2=r^4+S^2,
 $$
 where $a_2=l_a^{-2}+l_b^{-2}+l_c^{-2}$ and $a_3=l_a^{-1}l_b^{-1}l_c^{-1}$. By substituting $pr$ for $S$, where $p$ is the semiperimeter, we obtain that

 $$
 4a_2r^2p^2-8a_3r^3p^2=r^2+p^2.
 $$
 Since a general polynomial equation of degree four or lower can be solved in radicals \cite[page~126--130]{Postnikov}, the following statement holds.

 \begin{lemma}\label{L1}
{\it For every triangle $ABC$ with angle bisectors $l_a$, $l_b$, $l_c$,  all of the following objects are either expressible or not expressible by radicals of $l_a$, $l_b$, $l_c$ simultaneously:

\begin{itemize}
 \item the radius of the incircle of $ABC$,
 \item the area of $ABC$,
\item the perimeter of $ABC$.
\end{itemize}
}
 \end{lemma}
 In \cite{vRF}, it is also stated that the radius of the incircle of a triangle is a root of a polynomial of degree 16 with some rational functions of the angle bisectors for coefficients.

 After almost 100 years since the publication of \cite{vRF}, two papers were published in the same journal. We will see that a negative answer to the problem of expressing the area of a triangle in
radicals of the lengths of the angle bisectors  follows from these papers.

 The first of these papers was published in 1937 by  H. Wolff~\cite{Wolff}. It contains a proof, based on geometric reasons, that
 $\displaystyle\frac{1}{2r}$  is a root of the following polynomial:
 \begin{multline*}
  W(t)=t^{10}-\frac{5}{2}a_2t^8+\frac{7}{2}a_3t^7+\frac{33}{16}a_2^2t^6
  -\frac{47}{8}a_2a_3t^5+\\
  +\left(\frac{1}{4}a_2a_4-\frac{5}{8}a_2^3+\frac{61}{16}a_3^2\right)t^4+\left(\frac{5}{2}a_2^2a_3-\frac{1}{4}a_4a_3\right)t^3+\\
+\left(\frac{1}{16}a_2^4-\frac{1}{4}a_2^2a_4-\frac{25}{8}a_2a_3^2\right)t^2+\left(\frac{1}{2}a_2a_3a_4-\frac{1}{8}a_2^3a_3+\frac{5}{4}a_3^3\right)t+\\
+\left(\frac{1}{16}a_2^2a_3^2-\frac{1}{4}a_4a_3^2\right),
\end{multline*}
where
$$
a_2=l_a^{-2}+l_b^{-2}+l_c^{-2},\quad a_3=l_a^{-1}l_b^{-1}l_c^{-1},\quad \text{ and }\quad a_4=l_a^{-2}l_b^{-2}+l_b^{-2}l_c^{-2}+l_c^{-2}l_a^{-2}.
$$
It is also proved that the polynomial $W(t)$ is irreducible over the field $\mathbb{Q}(a_2,a_3,a_4)$. %Любопытно, что корнями многочлена $W(t)$ являются также
A year later B.L. van der Waerden \cite{vdWaerden} showed that the Galois group of the Wolff polynomial $W(t)$ over $\mathbb{Q}(a_2,a_3,a_4)$ is isomorphic to $S_{10}$, the symmetric group of degree
10 (which is non-solvable). So the roots of a polynomial $W(t)$ cannot be expressed in radicals of $a_2$, $a_3$, $a_4$ (see, for instance, \cite[pages 89--90]{Postnikov}). These papers provides no
explicit conclusion on solvability of the equation  $W(t)=0$ over $\mathbb{Q}(l_a,l_b,l_c)$ in radicals. However, it is quite easy to prove the following:

\begin{lemma}\label{T2}
%\noindent{\bf ТЕОРЕМА.}
{\it Roots of the polynomial $W(t)$ cannot be expressed in radicals of $l_a,$~$l_b,$~$l_c$.}
 \end{lemma}

 PROOF. Assume that a root of the polynomial $W(t)$ can be expressed in radicals of $l_a$, $l_b$, $l_c$. Let us show that then it can also be expressed in radicals of $a_2$, $a_3$, $a_4$, which
contradicts the results by van der Waerden and Wolff. To do that, it suffices to show that $l_a$, $l_b$, $l_c$ themselves can be expressed in radicals of $a_2,a_3,a_4$.

 Consider the following polynomials:
 $$
 U(t)= t^3-a_2t^2+a_4t-a_3^2 \quad\text{ and }\quad V(t)=U(t^2). $$

By Vieta's formulas, the roots of the polynomial $U(t)$ are equal to $l_a^{-2},l_b^{-2},l_c^{-2}$, and the roots of $V(t)$ are
$$
\pm\frac{1}{l_a},\quad\pm\frac{1}{l_b},\quad\pm\frac{1}{l_c}.
$$
Cardano's formula \cite[page~130]{Postnikov} implies that the roots of the cubic equation $U(t)=0$ and the bicubic equation $V(t)=0$ can be expressed in radicals of the coefficients $a_2$, $a_3^2$,
$a_4$. Hence $l_a$, $l_b$, $l_c$ can be expressed in radicals of $a_2$, $a_3$, $a_4$. This contradiction proves the lemma.

\bigskip

Now the following statement is a direct corollary of Lemma~\ref{T2}.

\begin{theorem}\label{C3}
{\it The radius of the incircle, as well as the area and the perimeter of a triangle cannot be expressed in radicals of the lengths of this triangle's angle bisectors.}
 \end{theorem}

Using a direct computer calculation in \cite{Magma}, one can show that if $l_a=1$, $l_b=2$ and $l_c=3$ then the Galois group of the polynomial $W(t)$ over $\mathbb{Q}$ is isomorphic to~$S_{10}$. So
the following statement holds.

\begin{theorem}\label{C4}

{\it The radius of the incircle, the area and the perimeter of a triangle with the length of the angle bisectors $1$, $2$ and $3$ cannot be expressed in radicals of rational numbers.}
 \end{theorem}

 \bigskip

We would like to express our gratitude to the team of translators from Novosibirsk State University: V.\,A.\,Afanas'ev, K.\,A.\,Kaushan, D.\,V.\,Lytkina, and T.\,R.\,Nasybullov who translated our paper into English.

\end{document}